\documentclass{elsart3-1}

\usepackage{amssymb}

\usepackage[english,francais]{babel}

\newtheorem{theorem}{Theorem}[section]

\newtheorem{e-proposition}[theorem]{Proposition}

\newtheorem{e-definition}[theorem]{Definition\rm}

\renewcommand{\theequation}{\arabic{equation}}
\def\N{\mathbb N}
\def\C{\mathbb C}

\def\R{\mathbb R}
\def\T{\mathbb T}

\def\Z{\mathbb Z}

\def\og{\leavevmode\raise.3ex\hbox{$\scriptscriptstyle\langle\!\langle$~}}
\def\fg{\leavevmode\raise.3ex\hbox{~$\!\scriptscriptstyle\,\rangle\!\rangle$}}

\journal{the Acad\'emie des sciences}
\begin{document}
\centerline{}
\begin{frontmatter}


\selectlanguage{english}
\title{Effective stability for slow time-dependent near-integrable Hamiltonians and application}


\selectlanguage{english}
\author[authorlabel1]{Abed Bounemoura},
\ead{abedbou@gmail.com}

\address[authorlabel1]{Institut des Hautes \'Etudes Scientifiques, 35 route de Chartres 91440 Bures-sur-Yvette (abed@ihes.fr)}


\medskip
\begin{center}
{\small Received *****; accepted after revision +++++\\
Presented by £££££}
\end{center}

\begin{abstract}
\selectlanguage{english}
The aim of this note is to prove a result of effective stability for a non-autonomous perturbation of an integrable Hamiltonian system, provided that the perturbation depends slowly on time. Then we use this result to clarify and extend a stability result of Giorgilli and Zehnder for a mechanical system with an arbitrary time-dependent potential.
{\it To cite this article: }

\vskip 0.5\baselineskip

\selectlanguage{francais}
\noindent{\bf R\'esum\'e} \vskip 0.5\baselineskip \noindent
{\bf Stabilit\'e effective pour des Hamiltoniens presque int\'egrables lentement non-autonomes et application. }
Le but de cette note est de d\'emontrer un r\'esultat de stabilit\'e effective pour une perturbation non-autonome d'un syst\`eme hamiltonien int\'egrable, sous la condition que la perturbation d\'epende lentement du temps. Nous utilisons ensuite ce r\'esultat pour clarifier et g\'en\'eraliser un r\'esultat de stabilit\'e de Giorgilli et Zehnder pour des syst\`emes m\'ecaniques dont le potentiel d\'epend arbitrairement du temps. 
{\it Pour citer cet article~: }

\end{abstract}
\end{frontmatter}

\selectlanguage{english}
\section{Introduction}\label{s1}

Let $n \in \N$, $n \geq 2$, $\T^n=\R^n / \Z^n$ and consider the Hamiltonian system defined by $H : \T^n \times \R^n \times \R \rightarrow \R$,
\begin{equation}\label{H}
H(\theta,I,t)=h(I)+\varepsilon f(\theta,I,t), \quad (\theta,I,t)=(\theta_1,\dots,\theta_n,I_1,\dots,I_n,t)\in \T^n \times \R^n \times \R, \quad \varepsilon>0.   
\end{equation}
Nekhoroshev proved (\cite{1}) that whenever $h$ is steep (see \S\ref{s2} for a definition), $f(\theta,I,t)=f(\theta,I)$ is time-independent and $H$ is real-analytic, there exist positive constants $\varepsilon_0,c_1,c_2,c_3,a,b$ such that for $\varepsilon \leq \varepsilon_0$, for all solutions $(\theta(t),I(t))$ we have the following stability estimates
\begin{equation}\label{estim}
|I(t)-I(0)|=\max_{1 \leq i \leq n}|I_i(t)-I_i(0)| \leq c_1\varepsilon^b, \quad |t|\leq c_2\exp(c_3\varepsilon^{-a}). 
\end{equation}
In the particular case where $h$ is (strictly uniformly) convex or quasi-convex, following a work of Lochak (\cite{4}) it was proved (\cite{2},\cite{3}), using preservation of energy arguments, that one can choose $a=b=(2n)^{-1}$ in~(\ref{estim}), and that these values are close to optimal. In the general steep case, however, there are still no realistic values for these stability exponents $a$ and $b$. The purpose of this note is to discuss to which extent stability estimates similar to~(\ref{estim}) hold true if the perturbation is allowed to depend on time.

Assume first that $f$ depends periodically on time, that is $f(\theta,I,t)=f(\theta,I,t+1)$ in~(\ref{H}). Removing the time-dependence by adding an extra degree of freedom, the Hamiltonian is equivalent to
\[ \tilde{H}(\theta,\varphi,I,J)=\tilde{h}(I,J)+\varepsilon f(\theta,\varphi,I), \quad (\theta,\varphi=t,I,J)\in \T^{n}\times\T \times \R^{n} \times \R, \quad \tilde{h}(I,J)=h(I)+J. \]
It turns out that if $h$ is convex, then $\tilde{h}$ is quasi-convex and so~(\ref{estim}) holds true with $a=b=(2(n+1))^{-1}$. In general, it is possible for $\tilde{h}$ to be steep in which case~(\ref{estim}) is satisfied, but it is not clear how to formulate a condition on $h$ (and not on $\tilde{h}$) to ensure that~(\ref{estim}) holds true.   

Now assume that $f$ depends quasi-periodically on time, that is $f(\theta,I,t)=f(\theta,I,t\omega)$ in~(\ref{H}) for some vector $\omega \in \R^m$ which we can assume to be non resonant ($k\cdot\omega \neq 0$ for any non-zero $k\in\Z^m$) and $f : \T^n \times \R^n \times \T^m \rightarrow \R$. As before, the time-dependence can be removed by adding $m$ degrees of freedom and we are led to consider $\tilde{H}(\theta,\varphi,I,J)=\tilde{h}(I,J)+\varepsilon f(\theta,\varphi,I)$ but this time
\[ (\theta,\varphi=t\omega,I,J)\in \T^{n} \times \T^m \times \R^{n}\times \R^{m}, \quad \tilde{h}(I,J)=h(I)+\omega\cdot J. \]
It was conjectured by Chirikov (\cite{5}), and then again by Lochak (\cite{6}), that if $h$ is convex and $\omega$ satisfies a Diophantine condition of exponent $\tau\geq m-1$ (there exists a constant $\gamma>0$ such that $|k.\omega|\geq \gamma|k|^{-\tau}$ for any non-zero $k\in\Z^m$), then the estimates~(\ref{estim}) hold true and moreover we can choose $a=b=(2(n+1+\tau))^{-1}$. If $m=1$, then $\tau=0$ and we are in the periodic case so the conjecture is true. However, if $m>1$, $\tilde{h}$ cannot be steep and the problem is still completely open. Even though the conjecture is sometimes considered as granted (for instance in \cite{7}), there is still no proof. Needless to say that the situation in the general case (without the convexity assumption on $h$) is even more complicated.

In a different direction, Giorgilli and Zehnder (\cite{8}) considered the following time-dependent Hamiltonian
\[ G(\theta,I,t)=h_2(I)+V(\theta,t), \quad (\theta,I,t)\in \T^{n} \times \R^{n} \times \R, \quad h_2(I)=I_1^2+\cdots+I_n^2, \]
and proved the following Nekhoroshev type result: if $G$ is real-analytic and $V$ is uniformly bounded, then for $R$ sufficiently large, if $I_0$ belongs to the ball $B_R$ of radius $R$ centered at the origin, then $I(t) \in B_{2R}$ for $|t|\leq c_2\exp(c_3 R^{d})$ for some positive constants $c_2,c_3$ and $d$. Even though such a system is clearly not of the form~(\ref{H}), the fact that no restriction on the time-dependence is imposed in their result has lead to several confusions. In \cite{8}, the authors themselves assert that ``extra work is needed because the time-dependence is not assumed to be periodic or quasi-periodic". Even more surprising, one can read (in \cite{9} for instance) that this result implies that the estimates~(\ref{estim}) hold true for~(\ref{H}) without any restriction on the time-dependence. Concerning the latter assertion, it is simply wrong and it seems very unlikely to have non trivial stability estimates for~(\ref{H}) with an arbitrary time-dependence. As for the former assertion, it is not difficult to see that the system considered in \cite{8} can be given the form~(\ref{H}), with a perturbation depending slowly on time, in their example $f(\theta,I,t)=f(\theta,I,\varepsilon^{1/2}t)$. We will show in \S\ref{s2} that for such Hamiltonian systems depending slowly on time, essentially classical techniques can be used to prove that~(\ref{estim}) hold true, and that the non-periodicity or non-quasi-periodicity of time plays absolutely no role (as a matter of fact, we already explained that for a periodic or quasi-periodic time dependence which is not slow, basic questions are still open). Then, in \S\ref{s3}, we will use this result to derive, in a simpler way, a more general statement than the one contained in \cite{8}.

\section{A stability result}\label{s2}

For a given $\rho>0$, a function $h \in C^2(B_\rho)$ is said to be steep if for any $k \in \{1, \dots, n\}$, there exist positive constants $p_k,C_k,\delta_k$ such that for any affine subspace $\lambda_k$ of dimension $k$ intersecting $B_\rho$, and any continuous curve $\gamma : [0,1] \rightarrow \lambda_k \cap B$ with $|\gamma(0)-\gamma(1)|=\delta<\delta_k$, there exists $t_*\in [0,1]$ such that $|\gamma (t)-\gamma (0)|< \delta$ for all $t\in [0,t_*]$ and
$|\Pi_{\Lambda_k}(\nabla h(\gamma (t_*)))| > C_k \delta^{p_k}$, where $\Pi_{\Lambda_k}$ is the projection onto $\Lambda_k$, the direction of $\lambda_k$. Then, given $r,s>0$, let us define the complex domain
\[ \mathcal{D}_{r,s}=\{(\theta,I,t)\in(\C^n/\Z^n)\times \C^{n} \times \C \; | \; |(\mathrm{Im}(\theta_1),\dots,\mathrm{Im}(\theta_n))|<s, \; |\mathrm{Im}(t)|<s,\; d(I,B_\rho)< r\},  \]
and for a fixed constant $1/2 \leq c \leq 1$, we consider $H(\theta,I,t)=h(I)+\varepsilon f(\theta,I,\varepsilon^ct)$ defined on $\mathcal{D}_{r,s}$ and real-analytic (that is $H$ is analytic and real-valued for real arguments). Finally, we assume that there exists a positive constant $M$ such that the operator norm $|\nabla^2 h(I)|\leq M$ for any $I\in B_\delta$, and that $|f(\theta,I,t)|\leq 1$ for any $(\theta,I,t) \in \mathcal{D}_{r,s}$.  

\medskip

\begin{theorem}\label{th1}
Under the previous assumptions, there exist positive constants $\varepsilon_0,c_1,c_2,c_3,$ that depend on $n,\rho,p_k,C_k,\delta_k,r,s,M$, and positive constants $a,b$ that depend only on $n,p_k$, such that if $\varepsilon \leq \varepsilon_0$, for all solutions $(\theta(t),I(t))$ of the Hamiltonian system defined by $H$, if $I(0)\in B_{\rho/2}$, then $|I(t)-I(0)| \leq c_1\varepsilon^b$ for all $|t|\leq c_2\exp(c_3\varepsilon^{-a})$. 
\end{theorem}

\medskip

Let us explain the proof, which follows from the arguments exposed in \cite{10} or \cite{11}, up to some technical points we shall detail now. First we remove the time-dependence: we let $x=\varepsilon^ct$ and we introduce a variable $y$ canonically conjugated to $x$, so that the Hamiltonian is equivalent to
\begin{equation}\label{HamD}
\tilde{H}(\theta,I,x,y)=h(I)+\varepsilon^c y +\varepsilon f(\theta,I,x)=h(I)+\tilde{f}(\theta,I,x,y), \quad (\theta,I,x,y) \in \tilde{\mathcal{D}}_{r,s}, 
\end{equation}
where $\tilde{\mathcal{D}}_{r,s}=\mathcal{D}_{r,s} \times\{y \in \C \; | \; |\mathrm{Im}(y)|<s \}$. The fact that the dependence on time is slow allows us to keep the integrable part fixed when removing the time-dependence, as one can consider that the extra degree of freedom only affects the perturbation. The new perturbation $\tilde{f}$ depends on parameters or ``degenerate" variables $x$ and $y$ (degenerate since they are not present in the integrable part), and such systems were already considered by Nekhoroshev (\cite{1}). However, a difficulty arise: for subsequent arguments, it is important for the (real part of the) variable $y$ to be unbounded, which is indeed the case by our definition of $\tilde{\mathcal{D}}_{r,s}$; but on this extended domain $\tilde{f}$ is unbounded and this implies that $\tilde{H}$ in~(\ref{HamD}) is not a perturbation of $h$. Yet the Hamiltonian vector field $X_{\tilde{H}}$ can be considered as a perturbation of $X_{h}$, as $X_{\tilde{f}}=(\partial_{I}\tilde{f},-\partial_{\theta}\tilde{f},\partial_{y}\tilde{f},-\partial_{x}\tilde{f})=(\varepsilon \partial_{I}f,-\varepsilon \partial_{\theta}f,\varepsilon^c,-\partial_{t}f)$, and so $X_{\tilde{f}}$ is bounded (by a Cauchy estimate) on the domain $\tilde{\mathcal{D}}_{r/2,s/2}$ for instance. As a consequence, even when $h$ is convex one cannot use preservation of energy arguments as it is the case in~\cite{2},~\cite{3},~\cite{4}, and in general one has use a perturbation theory that deals only with vector fields: the proofs in \cite{10} and \cite{11} accommodate both requirements. Now the analytic part of \cite{10} and \cite{11} goes exactly the same way for~(\ref{HamD}) by simply considering $x$ and $y$ as ``dummy" variables: given an integer parameter $m\geq 1$ which will be determined by the geometric part in terms of $\varepsilon$, on suitable domains resonant normal forms with a remainder of size bounded by a constant times $\varepsilon^ce^{-m}$ are constructed (note that the size of the perturbation $X_{\tilde{f}}$ is of order $\varepsilon^c$ and $c\leq 1$, but its ``effective" size is of order $\varepsilon$ and so $m$ will be determined in terms of $\varepsilon$ and not $\varepsilon^c$; $\varepsilon^c$ just enters the pre-factor in the exponential and will not alter the radius of confinement $\varepsilon^b$ as we always have $b\leq 1/2$ whereas $c\geq 1/2$.). The geometric part of \cite{10} and \cite{11} also goes exactly the same way since the time of escape (of the domain) of the degenerate variables $x$ and $y$ is infinite (as the domain is unbounded in these directions), $m$ is eventually chosen proportional to $\varepsilon^{-a}$, the radius of confinement is chosen proportional to $\varepsilon^b$ and the stability time is bounded by a constant times $e^{-m}$. 

Now let us add two remarks on the statement of Theorem~\ref{th1}. First, the exponents $a$ and $b$ are the same as in~(\ref{estim}) when the perturbation is time-independent. It is reasonable to expect that if $h$ is convex, then $a=b=(2n)^{-1}$ in Theorem~\ref{th1}, but we already explained that we cannot use preservation of energy arguments and so we cannot reach these values: the problem actually reduces to the problem of finding realistic values of $a$ and $b$ in the general steep case, which is still open. Then, using~\cite{10} and~\cite{11}, the statement of Theorem~\ref{th1} can be generalized in two ways: using~\cite{10} the statement holds true for the much wider class of Diophantine steep functions introduced by Niederman (which is a prevalent class of functions), using~\cite{11} the statement holds true for $\alpha$-Gevrey Hamiltonians for $\alpha\geq 1$ (with $\exp(c_3\varepsilon^{-a})$ replaced by $\exp(c_3\varepsilon^{-\alpha^{-1}a})$, recall that $1$-Gevrey is real-analytic) and for $C^k$ Hamiltonians, $k\geq n+1$ (with $\exp(c_3\varepsilon^{-a})$ replaced by $c_3\varepsilon^{-k^*a}$, with $k^*$ the largest integer $l\geq 1$ such that $k\geq ln+1$).

\section{An application}\label{s3}

Now we come back to the problem studied in \cite{8}, and more generally we consider, for an integer $p\geq 2$,
\[ G(\theta,I,t)=h_p(I)+V(\theta,t), \quad (\theta,I,t)\in \T^{n} \times \R^{n} \times \R, \quad h_p(I)=I_1^p+\cdots+I_n^p. \]
The case $p=2$ corresponds to~\cite{8} and $h_2$ is convex, for $p\geq 3$ the function $h_p$ is not convex but it is steep with $p_k=p-1$ and $C_k=\delta_k=1$ for all $1 \leq k \leq n$. The function $V$ is assumed to be real-analytic, defined on $\mathcal{D}_{s}=\{(\theta,t)\in(\C^n/\Z^n)\times \C \; | \; |(\mathrm{Im}(\theta_1),\dots,\mathrm{Im}(\theta_n)|<s, \; |\mathrm{Im}(t)|<s\}$, and it is assumed that $|V(\theta,t)|\leq 1$ for all $(\theta,t)\in \mathcal{D}_{s}$.

\medskip

\begin{theorem}\label{th2}
Under the previous assumptions, there exist positive constants $R_0,c_1,c_2,c_3$ that depend on $n,p,s$, and positive constants $a',b'$ that depend only on $n,p$, such that if $R\geq R_0$, for all solutions $(\theta(t),I(t))$ of the Hamiltonian system defined by $G$, if $I(0)\in B_{R}$, then $|I(t)-I(0)| \leq c_1R^{1-b'}$ for all $|t|\leq c_2R^{1-p}\exp(c_3R^{a'})$. 
\end{theorem}

\medskip

The proof is a direct application of Theorem~\ref{th1}. Indeed, for $R>0$ consider the scalings
\[ I=RI', \quad \theta=\theta', \quad G=R^p G', \quad t=R^{1-p}t'. \]
Then the Hamiltonian $G(\theta,I,t)$, for $(\theta,I,t) \in \mathcal{D}_{s} \times B_{2R}$, is equivalent to the Hamiltonian $G'(\theta',I',t')$, for $(\theta',I',t') \in \mathcal{D}_{s} \times B_{2}$, where $G'(\theta',I',t')=h_p(I')+R^{-p} V(\theta',R^{1-p}t')$. Hence we can apply Theorem~\ref{th1} to the Hamiltonian $G'$, with $\varepsilon=R^{-p}$, $c=(p-1)p^{-1}$, $\rho=2$ and $M$ which depends only on $p$: there exist positive constants $\varepsilon_0,c_1,c_2,c_3,$ that depend on $n,p,s$, and positive constants $a,b$ that depend only on $n,p$, such that if $\varepsilon \leq \varepsilon_0$, for all solutions $(\theta'(t'),I'(t'))$ of the Hamiltonian system defined by $G'$, if $I'(0)\in B_{1}$, then $|I'(t')-I'(0)| \leq c_1\varepsilon^b$ for all $|t'|\leq c_2\exp(c_3\varepsilon^{-a})$. Recalling that $\varepsilon=R^{-p}$, this means that if $R \geq R_0=\varepsilon_0^{-p^{-1}}$, for all $I'(0)\in B_{1}$ we have $|I'(t')-I'(0)| \leq c_1R^{-b'}$ for all $|t'|\leq c_2\exp(c_3R^{a'})$ for $b'=pb$ and $a'=pa$. Now scaling back to the original variables, for all $I(0)\in B_{R}$, we have $|I(t)-I(0)| \leq c_1R^{1-b'}$ for all $|t|\leq c_2R^{1-p}\exp(c_3R^{a'})$.

Now let us add some comments on the statement of Theorem~\ref{th2}. The estimate $|I(t)-I(0)| \leq c_1R^{1-b'}$ is stronger than $|I(t)-I(0)| \leq R$ (as in the argument above, the estimate $|I'(t')-I'(0)| \leq c_1\varepsilon^b$ is stronger than $|I'(t')-I'(0)| \leq 1$) and hence it is stronger than $I(t)\in B_{2R}$ if $I(0)\in B_R$. Moreover, we have $|t|\leq c_2\exp(c_3'R^{a'}) \leq c_2R^{1-p}\exp(c_3R^{a'})$ by restricting $c_3$ to a smaller value $c_3'$ and enlarging $R_0$ if necessary. So even for the convex case $p=2$ our statement is more accurate than the statement in~\cite{8}. In fact, for $p=2$, we already explained that we believe we can choose $a=b=(2n)^{-1}$, in which case the statement of Theorem~\ref{th2} would read $|I(t)-I(0)| \leq c_1R^{1-n^{-1}}$ for all $|t|\leq c_2R^{-1}\exp\left(c_3R^{n^{-1}}\right)$, which would be in perfect agreement with the much simpler autonomous case $V(\theta,t)=V(\theta)$ described in \cite{3}.        


\end{document}